\documentclass{article}[12pt]
\usepackage{amssymb}
\usepackage{eurosym}
\usepackage{epsfig}
\usepackage{amsmath}
\usepackage{amsfonts}
\usepackage{pgf,tikz}
\usepackage{enumerate}
\usepackage{cite}
\usepackage[colorlinks=true, linkcolor=black, citecolor=black, urlcolor=black]{hyperref}

\newcommand{\R}{{{\Bbb R}}}

\newtheorem{theorem}{\sc Theorem}[section]
\newtheorem{proposition}{\sc Proposition}[section]
\newtheorem{lemma}{\sc Lemma}[section]

\newtheorem{remark}{\sc Remark}[section]

\newtheorem{example}{\sc Example}[section]

\def\qed{\hbox to 0pt{}\hfill$\rlap{$\sqcap$}\sqcup$\medbreak}

\title{A fixed point index approach to Krasnosel'ski\u{\i}-Precup fixed point theorem in cones and applications}

\author{Jorge Rodr\'iguez--L\'opez} 

\date{}
\begin{document}
 \maketitle

\begin{center}  {\small CITMAga \& Departamento de Estat\'{\i}stica, An\'alise Matem\'atica e Optimizaci\'on, \\ Universidade de Santiago de Compostela, \\ 15782, Facultade de Matem\'aticas, Campus Vida, Santiago, Spain.\\  Email: jorgerodriguez.lopez@usc.es}
\end{center}

\medbreak

\noindent {\it Abstract.} We present an alternative approach to the vector version of Krasnosel'ski\u{\i} compression-expansion fixed point theorem due to Precup, which is based on the fixed point index. It allows us to obtain new general versions of this fixed point theorem and also multiplicity results. We emphasize that all of them are coexistence fixed point theorems for operator systems, that means that every component of the fixed points obtained is non-trivial. Finally, these coexistence fixed point theorems are applied to obtain results concerning the existence of positive solutions for systems of Hammerstein integral equations and radially symmetric solutions of  $(p_1,p_2)$-Laplacian systems.

\medbreak

\noindent     \textit{2020 MSC:} 47H10, 47H11, 45G15, 34B18, 35J92.

\medbreak

\noindent     \textit{Key words and phrases.} Coexistence fixed point; fixed point index; positive solution; Hammerstein systems; $p$-Laplacian system; radial solution.

\section{Introduction}

Krasnosel'ski\u{\i} compression-expansion fixed point theorem is between the main tools of Nonlinear Analysis for proving the existence of non-trivial solutions of different types of boundary value problems. Basically, assuming cone-compression or cone-expansion conditions on the boundary of an annulus, it ensures the existence of fixed points of compact operators defined in cones of normed linear spaces.

In the case of systems, the localization of the fixed points obtained by means of Krasnosel'ski\u{\i} theorem is not given independently in each component, so there is not guarantee that all the components of the fixed point are non-trivial, as already pointed out in \cite{CZ,PrecupFPT,PrecupSDC}. This fact motivated Precup to establish the \textit{vector version} of Krasnosel'ski\u{\i} fixed point theorem \cite{PrecupFPT,PrecupSDC} (see Theorem \ref{th_KP} below), which provides a component-wise localization of the fixed points. Thus, it gives sufficient conditions for the existence of a \textit{coexistence fixed point} as coined by Lan \cite{Lan1}, that is, a fixed point with all the components different from zero. As far as we know, there are only a few papers in the literature which deal with theoretical results concerning coexistence fixed points of compact maps, see \cite{CZ,Lan1,PrecupFPT,PrecupSDC,PreRod}. 

Moreover, under the assumptions of the vector version of Krasnosel'ski\u{\i} fixed point theorem, each component of the compact map may have a different behavior, namely, compression or expansion (see Remark \ref{rmk_posibCE} below). To the best of our knowledge, fixed point theorems for expansive-compressive maps are not common in the literature, we refer the interested reader to the results due to Mawhin \cite{Maw_CE} in this direction.

In this paper, we will refer to the mentioned vector version of Krasnosel'ski\u{\i} compression-expansion fixed point theorem due to Precup as Krasnosel'ski\u{\i}-Precup fixed point theorem in cones. It is well-known that the classical Krasnosel'ski\u{\i} fixed point theorem can be proved via fixed point index for compact maps, so our aim here is to present an alternative approach to Krasnosel'ski\u{\i}-Precup fixed point theorem based on this tool. It has its own interest since 
\begin{enumerate}
	\item[a)] it provides a way to extend Krasnosel'ski\u{\i}-Precup fixed point theorem to operators defined in more general domains;
	\item[b)] the computation of the fixed point index allows to obtain easily new multiplicity results;
	\item[c)] the proof can be replicated for other classes of maps for which a fixed point index theory is available (as, for instance, upper semicontinuous multivalued maps).
\end{enumerate} 
In addition, we prove that the finite-dimensional version of Krasnosel'ski\u{\i}-Precup fixed point theorem is equivalent to Poincar\'e-Miranda zeros theorem. This fact gives a connection between it and classical results. 

In the last section, we apply the theoretical results obtained in Section \ref{sec_KP} to two different problems: systems of Hammerstein integral equations and radially symmetric solutions of Dirichlet problems for $(p_1,p_2)$-Laplacian systems. Note that Krasnosel'ski\u{\i}-Precup fixed point theorem has been already employed by several authors in order to study the existence, localization and multiplicity of positive solutions for different types of systems of boundary value problems, see for instance \cite{feng,PrecupJMAA,WCS,WR,ZH}. Our intention is to emphasize the applicability of the new fixed point theorems established here and so we present a multiplicity result for a system of Hammerstein type equations, which complements previous results in the literature, see \cite{CCI14,fig_tojo,im,Lan1,PreRod} and the references therein. Moreover, concerning radial solutions of $(p_1,p_2)$-Laplacian systems, our sufficient conditions provide not only the existence of positive solutions, but also a novel localization of them, cf. \cite{OrW,ZH}. 

\section{Krasnosel'ski\u{\i}-Precup fixed point theorem in cones}\label{sec_KP}

First, we recall Krasnosel'ski\u{\i} compression-expansion fixed point theorem in cones. 

In the sequel, we need the following notions. A closed convex subset $K$ of a normed linear space $(X,\left\|\cdot\right\|)$ is a \textit{cone} if $\lambda\,u\in K$ for every $u\in K$ and for all $\lambda\geq 0$, and $K\cap (-K)=\{0\}$.
A cone $K$ induces the partial order in $X$ given by $u\preceq v$ if and only if $v-u\in K$. Moreover, we shall say that $u\prec v$ if $v-u\in K\setminus\{0\}$. 

The following notations will be useful: for given $r,R\in \mathbb{R}_+:=[0,\infty)$, $0<r<R$, we define
\[K_{r,R}:=\{u\in K:r<\left\|u\right\|<R \} \quad \text{ and } \quad \overline{K}_{r,R}:=\{u\in K:r\leq\left\|u\right\|\leq R \}. \]

\begin{theorem}[Krasnosel'ski\u{\i}]
	Let $(X,\left\|\cdot\right\|)$ be a normed linear space, $K$ a cone in $X$ and $r,R\in\R_+$, $0<r<R$.
	
	Consider a compact map $T:\overline{K}_{r,R}\rightarrow K$ satisfying one of the following conditions:
	\begin{enumerate}
		\item[$(a)$] $T(u)\nprec u$ if $\left\|u\right\|=r$ and $T(u)\nsucc u$ if $\left\|u\right\|=R$;
		\item[$(b)$] $T(u)\nsucc u$ if $\left\|u\right\|=r$ and $T(u)\nprec u$ if $\left\|u\right\|=R$.
	\end{enumerate}
	Then $T$ has at least a fixed point $u\in K$ with $r\leq\left\|u\right\|\leq R$.
\end{theorem}

Condition $(a)$ in Krasnosel'ski\u{\i} theorem is usually called a \textit{compression type condition}, whereas authors often refer to condition $(b)$ as the \textit{cone-expansion condition}. Similarly, in case $(a)$ we will say that the operator $T$ is \textit{compressive}, while in case $(b)$ $T$ is called an \textit{expansive} operator.

It is well-known that conditions $(a)$ and $(b)$ can be weakened as \textit{homotopy type conditions}. In this way, we have the \textit{homotopy version of Krasnosel'ski\u{\i} theorem} or \textit{Krasnosel'ski\u{\i}-Benjamin theorem}, see for instance \cite{amann}.

\begin{theorem}[Krasnosel'ski\u{\i}-Benjamin]
	Let $(X,\left\|\cdot\right\|)$ be a normed linear space, $K$ a cone in $X$ and $r,R\in\R_+$, $0<r<R$.
	
	Assume that $T:\overline{K}_{r,R}\rightarrow K$ is a compact map and there exists $h\in K\setminus\{0\}$ such that one of the following conditions is satisfied:
	\begin{enumerate}
		\item[$(a)$] $T(u)+\mu\,h\neq u$ if $\left\|u\right\|=r$ and $\mu>0$, and $T(u)\neq \lambda\, u$ if $\left\|u\right\|=R$ and $\lambda>1$;
		\item[$(b)$] $T(u)\neq \lambda\, u$ if $\left\|u\right\|=r$ and $\lambda>1$ , and $T(u)+\mu\,h\neq u$ if $\left\|u\right\|=R$ and $\mu>0$. 
	\end{enumerate}
	Then $T$ has at least a fixed point $u\in K$ with $r\leq\left\|u\right\|\leq R$.
\end{theorem}

In \cite{PrecupFPT,PrecupSDC}, Precup proposed a compression-expansion type fixed point theorem for systems of operators. The main novelty is that compression-expansion conditions are given in a component-wise manner, in what was called the \textit{vector version of Krasnosel'ski\u{\i} fixed point theorem}. Let us recall this result.

Consider two cones $K_1$ and $K_2$ of a normed linear space $X$, and so $K:=K_{1}\times K_{2}$ is a cone of $X^{2}=X\times X$. For $r,R\in \mathbb{R}_{+}^{2}$, $r=(r_{1},r_{2})$, $R=(R_{1},R_{2})$, with $0<r_i<R_i$ $(i=1,2)$, we denote
\begin{align*}
(\overline{K}_{i})_{r_{i},R_{i}}& :=\{u\in K_{i}\,:\,r_{i}\leq \left\Vert u\right\Vert
\leq R_{i}\}\quad (i=1,2), \\
\overline{K}_{r,R}& :=\{u=\left( u_{1},u_{2}\right) \in K\,:\,r_{i}\leq \left\Vert
u_{i}\right\Vert \leq R_{i}\quad \text{for }i=1,2\}.
\end{align*}
Clearly, $\overline{K}_{r,R}=(\overline{K}_{1})_{r_{1},R_{1}}\times (\overline{K}_{2})_{r_{2},R_{2}}$.

The aim of the vector version of Krasnosel'ski\u{\i} theorem is to obtain a solution $u=(u_1,u_2)$ to the operator system 
\[\left\{\begin{array}{l} u_1=T_1(u_1,u_2), \\ u_2=T_2(u_1,u_2), \end{array} \right. \]
located in the set $\overline{K}_{r,R}$, that is, $u=(u_1,u_2)\in K$ and $r_i\leq\left\|u_i\right\|\leq R_i$, $i=1,2$.

\begin{theorem}[Krasnosel'ski\u{\i}-Precup]\label{th_KP}
	Let $(X,\left\|\cdot\right\|)$ be a normed linear space, $K_1$ and $K_2$ two cones in $X$ and $r,R\in\R^2_+$, $r=(r_{1},r_{2})$, $R=(R_{1},R_{2})$, with $0<r_i<R_i$ $(i=1,2)$.
	
	Assume that $T=(T_1,T_2):\overline{K}_{r,R}\rightarrow K$ is a compact map and for each $i\in\{1,2\}$ there exists $h_i\in K_i\setminus\{0\}$ such that one of the following conditions is satisfied in $\overline{K}_{r,R}$:
	\begin{enumerate}
		\item[$(a)$] $T_i(u)+\mu\,h_i\neq u_i$ if $\left\|u_i\right\|=r_i$ and $\mu>0$, and $T_i(u)\neq \lambda\, u_i$ if $\left\|u_i\right\|=R_i$ and $\lambda>1$;
		\item[$(b)$] $T_i(u)\neq \lambda\, u_i$ if $\left\|u_i\right\|=r_i$ and $\lambda>1$, and $T_i(u)+\mu\,h_i\neq u_i$ if $\left\|u_i\right\|=R_i$ and $\mu>0$. 
	\end{enumerate}
	Then $T$ has at least a fixed point $u=(u_1,u_2)\in K$ with $r_i\leq\left\|u_i\right\|\leq R_i$ $(i=1,2)$.
\end{theorem}

\begin{remark}\label{rmk_posibCE}
	As already pointed out in \cite{PrecupFPT,PrecupSDC}, the operator $T$ may exhibit a different behavior (compression or expansion) in each component. More exactly, the following options are possible:
	\begin{enumerate}
		\item[$(i)$] both operators $T_1$ and $T_2$ are compressive;
		\item[$(ii)$] both operators $T_1$ and $T_2$ are expansive;
		\item[$(iii)$] one of the operators $T_1$ or $T_2$ is compressive, while the other one is expansive.
	\end{enumerate}
\end{remark}

Let us recall briefly the main ideas of the proof of Theorem \ref{th_KP} given in \cite{PrecupFPT,PrecupSDC}, which is essentially divided into two cases: 1) both operators $T_1$ and $T_2$ are compressive; 2) one of the operators is expansive. The first case relies on Schauder fixed point theorem. In the second case, the fixed point problem is reduced to an equivalent one in which both operators satisfy the compression type condition  and so the existence of a fixed point is ensured by the former case. 

Note that the same proof due to Precup remains valid for the following $n$-dimensional vector version of Krasnosel'ski\u{\i} fixed point theorem for an operator $T=(T_1,T_2,\dots,T_n)$ defined in $X^n$.

\begin{theorem}\label{th_KPn}
	Let $(X,\left\|\cdot\right\|)$ be a normed linear space, $K_1,\dots,K_n$ cones in $X$, $K:=K_1\times\cdots\times K_n$, $r,R\in\R^n_+$, $r=(r_{1},\dots,r_{n})$, $R=(R_{1},\dots,R_{n})$, with $0<r_i<R_i$ $(i=1,\dots,n)$, and $\overline{K}_{r,R} :=\{u=\left( u_{1},\dots,u_{n}\right) \in K\,:\,r_{i}\leq \left\Vert
	u_{i}\right\Vert \leq R_{i}\quad \text{for }i=1,\dots,n\}$.
	
	Assume that $T=(T_1,\dots,T_n):\overline{K}_{r,R}\rightarrow K$ is a compact map and for each $i\in\{1,\dots,n\}$ there exists $h_i\in K_i\setminus\{0\}$ such that one of the following conditions is satisfied in $\overline{K}_{r,R}$:
	\begin{enumerate}
		\item[$(a)$] $T_i(u)+\mu\,h_i\neq u_i$ if $\left\|u_i\right\|=r_i$ and $\mu>0$, and $T_i(u)\neq \lambda\, u_i$ if $\left\|u_i\right\|=R_i$ and $\lambda>1$;
		\item[$(b)$] $T_i(u)\neq \lambda\, u_i$ if $\left\|u_i\right\|=r_i$ and $\lambda>1$, and $T_i(u)+\mu\,h_i\neq u_i$ if $\left\|u_i\right\|=R_i$ and $\mu>0$. 
	\end{enumerate}
	Then $T$ has at least a fixed point $u=(u_1,\dots,u_n)\in K$ with $r_i\leq\left\|u_i\right\|\leq R_i$ $(i=1,\dots,n)$.
\end{theorem}

We highlight that our proof here is completely different to that due to Precup, since it is based on fixed point index theory independently of the possibility $(i)$--$(iii)$ in Remark \ref{rmk_posibCE}. In particular, our approach does not require to turn the expansive operators into compressive ones. 

\subsection{Fixed point index computation}

First, let us recall some of the useful properties of the fixed point index for compact maps. For more details, we refer the reader to \cite{amann,Deim,GraDug} (see also \cite{Inf}).

\begin{proposition}\label{prop_index}
	Let $P$ be a cone of a normed linear space, $U\subset P$ be a bounded relatively open set and $T:\overline{U}\rightarrow P$ be a compact map such that $T$ has no fixed points on the boundary of $U$ (denoted by $\partial\, U$). Then the fixed point index of $T$ in $P$ over $U$, $i_{P}(T,U)$, has the following properties:
	\begin{enumerate}
		\item (Additivity) Let $U$ be the disjoint union of two open sets $U_1$ and $U_2$. If $0\not\in(I-T)(\overline{U}\setminus(U_1\cup U_2))$, then \[i_{P}(T,U)=i_{P}(T,U_1)+i_{P}(T,U_2).\]
		\item (Existence) If $i_{P}(T,U)\neq 0$, then there exists $u\in U$ such that $u=Tu$.
		\item (Homotopy invariance) If $H:\overline{U}\times[0,1]\rightarrow P$ is a compact homotopy and $0\not\in(I-H)(\partial\,U\times[0,1])$, then
		\[i_{P}(H(\cdot,0),U)=i_{P}(H(\cdot,1),U).\]
		\item (Normalization) If $T$ is a constant map with $T(u)=u_0$ for every $u\in\overline{U}$, then
		\[i_{P}(T,U)=\left\{\begin{array}{ll} 1, & \text{ if } u_{0}\in U, \\ 0, & \text{ if } u_{0}\not\in\overline{U}. \end{array} \right. \]
	\end{enumerate}
\end{proposition}

Moreover, we have the following conditions concerning the computation of the fixed point index.

\begin{proposition}\label{prop_ind01}
	Let $U$ be a bounded relatively open subset of a cone $P$ such that $0\in U$ and $T:\overline{U}\rightarrow P$ be a compact map.
	\begin{enumerate}[$(a)$]
		\item If $T(u)\neq \lambda\,u$ for all $u\in \partial\, U$ and all $\lambda\geq 1$, then $i_{P}(T,U)=1$.
		\item If there exists $h\in P\setminus\{0\}$ such that $T(u)+\lambda\, h\neq u$ for every $\lambda\geq 0$ and all $u\in \partial\, U$, then $i_{P}(T,U)=0$.
	\end{enumerate}	
\end{proposition}

\medskip

In the sequel, let $(X,\left\|\cdot\right\|_X)$ and $(Y,\left\|\cdot\right\|_Y)$ be normed linear spaces, $K_{1}\subset X$, $K_2\subset Y$ two cones and $K:=K_{1}\times K_{2}$ the corresponding cone of $X\times Y$. When no confusion may occur, both norms $\left\|\cdot\right\|_X$ and $\left\|\cdot\right\|_Y$ will be simply denoted by $\left\|\cdot\right\|$. 

Now, we present a technical result which will be crucial in the computation of the fixed point index in the main results of this section. It was already proven in \cite[Lemma 2.3]{PreRod}, but we include here the proof again for the reader's convenience. 

\begin{lemma}\label{lem_index0}
	Let $U$ and $V$ be bounded relatively open subsets of $K_1$ and $K_2$, respectively, such that $0\in U$. 
	
	Assume that  $T:\overline{U\times V}\rightarrow K$, $T=(T_1,T_2)$, is a compact map and there exists $h\in K_{2}\setminus\{0\}$ such that
	\begin{align}
	T_{1}(u,v)&\neq \lambda\,u \quad \text{for } u\in\partial_{K_1}U, \ v\in \overline{V} \text{ and } \lambda\geq 1; \label{eq1:lem}  \\
	T_{2}(u,v)+\mu\, h&\neq v \quad \text{for } u\in\overline{U}, \ v\in\partial_{K_2}V \text{ and } \mu\geq 0. \label{eq2:lem} 
	\end{align}
	Then $i_{K}(T,U\times V)=0$.
\end{lemma}

\noindent
{\bf Proof.} Consider the homotopy $H:\overline{U\times V}\times[0,1]\rightarrow K$ given by
\[H((u,v),t)=(t\,T_1(u,v),T_2(u,v)+(1-t)\mu_0\, h), \]
with $\mu_0>0$ big enough such that $v\neq T_2(0,v)+\mu_0 h$ for all $v\in \overline{V}$. Note that the existence of such a positive number $\mu_0$ is guaranteed since $\overline{V}$ is bounded and $T$ is compact.

Assumptions \eqref{eq1:lem} and \eqref{eq2:lem} guarantee that the homotopy function $H$ has no fixed points on $\partial_{K}(U\times V)$. Therefore, by the homotopy invariance property of the fixed point index, we have that 
\begin{equation}\label{eq_lem_index0}
i_{K}(H(\cdot,0),U\times V)=i_{K}(H(\cdot,1),U\times V)=i_{K}(T,U\times V).
\end{equation}
On the other hand, for $t=0$, the map $H((u,v),0)=(0,T_2(u,v)+\mu_0\, h)$ has no fixed points in $U\times V$. Indeed, if $(u,v)\in U\times V$ is such a fixed point, then $u=0$ and $v=T_2(0,v)+\mu_0 h$, a contradiction with the hypothesis about $\mu_0$.
Hence, $i_{K}(H(\cdot,0),U\times V)=0$ and so the conclusion follows from \eqref{eq_lem_index0}.
\qed

\begin{remark}
	Obviously, the roles that play $T_1$ and $T_2$ in the statement of Lemma \ref{lem_index0}, given by assumptions \eqref{eq1:lem} and \eqref{eq2:lem}, are interchangeable.
\end{remark}

For $r,R\in\mathbb{R}^2_+$, $0<r_i<R_i$ ($i=1,2$), fixed, our aim is to compute the fixed point index of a compact operator $T=(T_1,T_2):\overline{K}_{r,R}\rightarrow K$ in the relatively open set
\[K_{r,R}:=\left\{u=\left( u_{1},u_{2}\right) \in K\,:\,r_{i}< \left\Vert
u_{i}\right\Vert < R_{i}\quad \text{for }i=1,2 \right\} \]
under the conditions of Krasnosel'ski\u{\i}-Precup fixed point theorem. Obviously, we need to assume also that $T$ has no fixed points on the boundary of $K_{r,R}$ in order to have the fixed point index well-defined over this set.

In the sequel, we will also use the following notations:
\begin{align*}
	(K_i)_{r_i}&=\{u\in K_{i}\,:\, \left\Vert u\right\Vert
	< r_{i}\}, \\ 
	(\overline{K}_i)_{r_i}&=\{u\in K_{i}\,:\, \left\Vert u\right\Vert
	\leq r_{i}\} \quad (i=1,2), \\
	K_{r}&=\left\{u=\left( u_{1},u_{2}\right) \in K\,:\left\Vert
	u_{i}\right\Vert < r_{i}\quad \text{for }i=1,2\right\}, \\ \overline{K}_{r}&=\{u=\left( u_{1},u_{2}\right) \in K\,:\left\Vert
	u_{i}\right\Vert \leq r_{i}\quad \text{for }i=1,2\}.
\end{align*}

To compute the fixed point index over $K_{r,R}$, we need to extend the definition of $T$ to the set $\overline{K}_{R}$. A key ingredient in our purpose is the definition of a retraction from $\overline{K}_{R}$ into $\overline{K}_{r,R}$. To do so, we use that $(\overline{K}_i)_{r_i,R_i}$ is a retract of $(\overline{K}_i)_{R_i}$ ($i=1,2$), see \cite[Example 3]{fel}. Indeed, we have the retraction $\rho_i: (\overline{K}_i)_{R_i}\rightarrow (\overline{K}_i)_{r_i,R_i}$ defined as
\begin{equation}\label{retract}
\rho_i(u_i)=\left\{\begin{array}{ll} r_i\displaystyle\frac{u_i+(r_i-\left\|u_i\right\|)^2 h_i}{\left\|u_i+(r_i-\left\|u_i\right\|)^2 h_i \right\|}, & \quad \text{if } \left\|u_i\right\|<r_i, \\[0.2cm] u_i,  & \quad \text{if } r_i\leq\left\|u_i\right\|\leq R_i,  \end{array} \right.
\end{equation}
where $h_i\in K_i\setminus\{0\}$ is fixed. Note that $\rho_i$ is well-defined: $\left\|u_i+(r_i-\left\|u_i\right\|)^2 h_i \right\|\neq 0$ for all $u_i\in (\overline{K}_i)_{r_i}$. Otherwise, $-u_i=(r_i-\left\|u_i\right\|)^2 h_i\in K_i$, what together with $u_i\in K_i$ implies $u_i=0$, from the definition of cone. Taking $u_i=0$, we have $\left\|r_i^2\,h_i\right\|>0$ since $r_i>0$ and $h_i\in K_i\setminus\{0\}$. Moreover, it is clear that $\rho_i$ is continuous and $\rho_i(u_i)=u_i$ for all $u_i\in (\overline{K}_i)_{r_i,R_i}$.

Now we are in a position to compute the fixed point index over the set $K_{r,R}$. First, we study the case in which both operators are compressive.

\begin{theorem}\label{th_CC}
%
	Assume that $T=(T_1,T_2):\overline{K}_{r,R}\rightarrow K$ is a compact map and for each $i\in\{1,2\}$ there exists $h_i\in K_i\setminus\{0\}$ such that the following conditions are satisfied in $\overline{K}_{r,R}$:
	\begin{enumerate}
		\item[$(i)$] $T_i(u)+\mu\,h_i\neq u_i$ if $\left\|u_i\right\|=r_i$ and $\mu\geq0$; 
		\item[$(ii)$] $T_i(u)\neq \lambda\, u_i$ if $\left\|u_i\right\|=R_i$ and $\lambda\geq 1$.
	\end{enumerate}
	Then \[i_{K}(T,K_{r,R})=1.\]
\end{theorem}

\noindent
{\bf Proof.} Consider the retraction $\rho:\overline{K}_{R}\rightarrow \overline{K}_{r,R}$ defined as 
$\rho(u_1,u_2)=\left(\rho_1(u_1),\rho_2(u_2)\right), $ where the functions $\rho_i$ are given by \eqref{retract}, $i=1,2$. Now, define the auxiliary map $N=(N_1,N_2):\overline{K}_{R}\rightarrow K$ as follows
\begin{equation}\label{eq_N}
N(u):=(T\circ\rho)(u). 
\end{equation}

Clearly, $N$ is a compact operator, $N(u)=T(u)$ for every $u\in\overline{K}_{r,R}$ and for each $i\in\{1,2\}$ the following conditions hold in $\overline{K}_{R}$:
\begin{enumerate}
	\item[$(i^*)$] $N_i(u)+\mu\,h_i\neq u_i$ if $\left\|u_i\right\|=r_i$ and $\mu\geq0$; 
	\item[$(ii^*)$] $N_i(u)\neq \lambda\, u_i$ if $\left\|u_i\right\|=R_i$ and $\lambda\geq 1$.
\end{enumerate}

Note that $(i^*)$ implies that $N(u)+\mu\,h\neq u$ for all $u\in\partial K_{r}$ and $\mu\geq 0$ (where $h=(h_1,h_2)$), so Proposition \ref{prop_ind01} yields \[i_{K}(N,K_r)=0. \]
Similarly, condition $(ii^*)$ in conjunction with Proposition \ref{prop_ind01} guarantee that \[i_{K}(N,K_R)=1. \]
Moreover, by Lemma \ref{lem_index0}, \[i_{K}(N,(K_1)_{R_1}\times(K_2)_{r_2})=i_{K}(N,(K_1)_{r_1}\times(K_2)_{R_2})=0. \]
By the additivity property of the fixed point index,
\[i_{K}(N,(K_1)_{r_1,R_1}\times(K_2)_{r_2})=i_{K}(N,(K_1)_{R_1}\times(K_2)_{r_2})-i_{K}(N,K_r)=0, \]
so we obtain
\[i_{K}(N,K_{r,R})=i_{K}(N,K_R)-i_{K}(N,(K_1)_{r_1,R_1}\times(K_2)_{r_2})-i_{K}(N,(K_1)_{r_1}\times(K_2)_{R_2})=1. \]

Finally, since $T=N$ on the set $\overline{K}_{r,R}$, we have $i_{K}(T,K_{r,R})=i_{K}(N,K_{r,R})=1$.
\qed

Next, let us consider the case in which we have compression for one operator and expansion for the other one. 

\begin{theorem}\label{th_CE}
%
	Assume that $T=(T_1,T_2):\overline{K}_{r,R}\rightarrow K$ is a compact map and for each $i\in\{1,2\}$ there exists $h_i\in K_i\setminus\{0\}$ such that the following conditions are satisfied in $\overline{K}_{r,R}$:
	\begin{enumerate}
		\item[$(i)$] $T_1(u)+\mu\,h_1\neq u_1$ if $\left\|u_1\right\|=r_1$ and $\mu\geq0$, and $T_1(u)\neq \lambda\, u_1$ if $\left\|u_1\right\|=R_1$ and $\lambda\geq 1$; 
		\item[$(ii)$] $T_2(u)+\mu\,h_2\neq u_2$ if $\left\|u_2\right\|=R_2$ and $\mu\geq0$, and $T_2(u)\neq \lambda\, u_2$ if $\left\|u_2\right\|=r_2$ and $\lambda\geq 1$.
	\end{enumerate}
	Then \[i_{K}(T,K_{r,R})=-1.\]
\end{theorem}

\noindent
{\bf Proof.} Consider the map $N:\overline{K}_{R}\rightarrow K$ defined as in \eqref{eq_N}. By Proposition \ref{prop_ind01}, 
\[i_{K}(N,(K_1)_{R_1}\times(K_2)_{r_2})=1, \qquad i_{K}(N,(K_1)_{r_1}\times(K_2)_{R_2} )=0. \]
Moreover, Lemma \ref{lem_index0} yields
\[i_{K}(N,K_r)=i_{K}(N,K_R)=0. \]

Hence, it follows from the additivity property of the index that
\[i_{K}(N,(K_1)_{r_1,R_1}\times(K_2)_{r_2})=i_{K}(N,(K_1)_{R_1}\times(K_2)_{r_2})-i_{K}(N,K_r)=1, \]
and so
\[i_{K}(N,K_{r,R})=i_{K}(N,K_R)-i_{K}(N,(K_1)_{r_1}\times(K_2)_{R_2} )-i_{K}(N,(K_1)_{r_1,R_1}\times(K_2)_{r_2})=-1. \]
Then, $i_{K}(T,K_{r,R})=-1$.
\qed

\begin{remark}
	The statement of Theorem \ref{th_CE} corresponds to the case in which $T_1$ is compressive and $T_2$ is expansive. Clearly, the same conclusion can be obtained if $T_1$ is expansive and $T_2$ is compressive.	
\end{remark}

Finally, we deal with the case in which both $T_1$ and $T_2$ are expansive.

\begin{theorem}\label{th_EE}
%
	Assume that $T=(T_1,T_2):\overline{K}_{r,R}\rightarrow K$ is a compact map and for each $i\in\{1,2\}$ there exists $h_i\in K_i\setminus\{0\}$ such that the following conditions are satisfied in $\overline{K}_{r,R}$:
	\begin{enumerate}
		\item[$(i)$] $T_i(u)+\mu\,h_i\neq u_i$ if $\left\|u_i\right\|=R_i$ and $\mu\geq0$; 
		\item[$(ii)$] $T_i(u)\neq \lambda\, u_i$ if $\left\|u_i\right\|=r_i$ and $\lambda\geq 1$.
	\end{enumerate}
	Then \[i_{K}(T,K_{r,R})=1.\]
\end{theorem}

\noindent
{\bf Proof.} Consider again the map $N:\overline{K}_{R}\rightarrow K$ defined as in \eqref{eq_N}. Now, by Proposition \ref{prop_ind01},
\[i_{K}(N,K_R)=0, \qquad i_{K}(N,K_r)=1. \]
Moreover, according with Lemma \ref{lem_index0}, 
\[i_{K}(N,(K_1)_{r_1}\times(K_2)_{R_2} )=i_{K}(N,(K_1)_{R_1}\times(K_2)_{r_2} )=0. \]

As a consequence of the additivity property of the index, we deduce that
\[i_{K}(N,(K_1)_{r_1}\times(K_2)_{r_2,R_2} )=i_{K}(N,(K_1)_{r_1}\times(K_2)_{R_2} )-i_{K}(N,K_r )=-1. \]
Therefore,
\[i_{K}(N,K_{r,R})=i_{K}(N,K_R)-i_{K}(N,(K_1)_{r_1}\times(K_2)_{r_2,R_2} )-i_{K}(N,(K_1)_{R_1}\times(K_2)_{r_2} )=1. \]
In conclusion, $i_{K}(T,K_{r,R})=1$, as wished.
\qed

\begin{remark}
	The computation of the fixed point index given by Theorems \ref{th_CC}, \ref{th_CE} and \ref{th_EE} is independent of that obtained in \cite[Theorem 2.8]{PreRod}. Indeed, in \cite{PreRod}, the assumptions on $T$ and, in particular, the homotopy conditions $(i)$--$(ii)$ were imposed in the whole set $\overline{K}_{R}$ instead of its subset $\overline{K}_{r,R}$, as here.
	
\end{remark}

As a straightforward consequence of the computation of the fixed point index provided by Theorems \ref{th_CC}, \ref{th_CE} and \ref{th_EE}, we have an alternative version of Theorem \ref{th_KP}.

\begin{theorem}\label{th_KP_index}
	Assume that $T=(T_1,T_2):\overline{K}_{r,R}\rightarrow K$ is a compact map and for each $i\in\{1,2\}$ there exists $h_i\in K_i\setminus\{0\}$ such that one of the following conditions is satisfied in $\overline{K}_{r,R}$:
	\begin{enumerate}
		\item[$(a)$] $T_i(u)+\mu\,h_i\neq u_i$ if $\left\|u_i\right\|=r_i$ and $\mu\geq 0$, and $T_i(u)\neq \lambda\, u_i$ if $\left\|u_i\right\|=R_i$ and $\lambda\geq 1$;
		\item[$(b)$] $T_i(u)\neq \lambda\, u_i$ if $\left\|u_i\right\|=r_i$ and $\lambda\geq 1$, and $T_i(u)+\mu\,h_i\neq u_i$ if $\left\|u_i\right\|=R_i$ and $\mu\geq 0$. 
	\end{enumerate}
	Then $T$ has at least a fixed point $u=(u_1,u_2)\in K$ with $r_i<\left\|u_i\right\|< R_i$ $(i=1,2)$.
\end{theorem}

\noindent
{\bf Proof.} By Theorems \ref{th_CC}, \ref{th_CE} and \ref{th_EE}, we have that
\[i_{K}(T,K_{r,R})=\pm 1, \]
and so the existence property of the fixed point index ensures that $T$ has at least a fixed point in $K_{r,R}$.
\qed

%

\subsection{Other versions of Krasnosel'ski\u{\i}-Precup fixed point theorem: different domains}

After the previous computation of the fixed point index of $T$ over the set $K_{r,R}$, we can think of proving similar results for operators $T$ defined in other regions different from $\overline{K}_{r,R}$. In this way, we increase the range of applicability of the original Krasnosel'ski\u{\i}-Precup fixed point theorem.

For each $i\in\{1,2\}$, let $\varphi_i:K_i\rightarrow\mathbb{R}_+$ be a continuous concave functional on $K_i$, that is, $\varphi_i$ is a continuous function and
\[\varphi_i\left(\lambda\,u+(1-\lambda)v\right)\geq \lambda\,\varphi_i(u)+(1-\lambda)\varphi_i(v), \quad \text{for all } u,v\in K_i, \ \lambda\in[0,1]. \]
Then, for $r,R\in\mathbb{R}^2_+$, $0<r_i<R_i$ ($i=1,2$), fixed, consider the sets
\begin{align*}
K^{\;\varphi}_{r,R}&:=\left\{u=(u_1,u_2)\in K\,:\,r_i< \varphi_i(u_i) \ \text{ and } \ \left\|u_i\right\|< R_i \quad \text{for }i=1,2 \right\}, \\ \overline{K}^{\;\varphi}_{r,R}&:=\left\{u=(u_1,u_2)\in K\,:\,r_i\leq \varphi_i(u_i) \ \text{ and } \ \left\|u_i\right\|\leq R_i \quad \text{for }i=1,2 \right\}.
\end{align*}
This type of sets has been already considered by Leggett and Williams \cite{LW} in the context of their celebrated fixed point theorem. Note that $\overline{K}^{\;\varphi}_{r,R}$ is a closed convex set. Hence, by Dugundji extension theorem (see \cite[Theorem 4.1]{dug} or \cite{GraDug}), the subset $\overline{K}^{\;\varphi}_{r,R}$ is a retract of $K$.

With this in mind, one can easily establish alternative versions of Theorems \ref{th_CC}, \ref{th_CE} and \ref{th_EE}. We sum up them in the following result.

\begin{theorem}\label{th_ind_funct}
%
	Assume that there exist continuous concave functionals $\varphi_i:K_i\rightarrow\mathbb{R}_+$ such that $\varphi_i(u)\leq\left\|u\right\|$ for all $u\in K_i$ ($i=1,2$), the set $K^{\;\varphi}_{r,R}$ is nonempty, $T=(T_1,T_2):\overline{K}^{\;\varphi}_{r,R}\rightarrow K$ is a compact map and for each $i\in\{1,2\}$ there exists $h_i\in K_i\setminus\{0\}$ such that one of the following conditions is satisfied in $\overline{K}^{\;\varphi}_{r,R}$:
	\begin{enumerate}
		\item[$(i)$] $T_i(u)+\mu\,h_i\neq u_i$ if $\varphi_i(u_i)=r_i$ and $\mu\geq 0$, and $T_i(u)\neq \lambda\, u_i$ if $\left\|u_i\right\|=R_i$ and $\lambda\geq 1$;
		\item[$(ii)$] $T_i(u)\neq \lambda\, u_i$ if $\varphi_i(u_i)=r_i$ and $\lambda\geq 1$, and $T_i(u)+\mu\,h_i\neq u_i$ if $\left\|u_i\right\|=R_i$ and $\mu\geq 0$. 
	\end{enumerate}
	
	Then \[i_{K}(T,K^{\;\varphi}_{r,R})=(-1)^{k},\]
	where $k=0$ if both $T_1$ and $T_2$ are compressive, $k=1$ if one of the operators $T_1$ or $T_2$ is compressive and the other one is expansive and $k=2$ if both operators are expansive.
\end{theorem}

\noindent
{\bf Proof.} Consider a retraction $\rho:K\rightarrow \overline{K}^{\;\varphi}_{r,R}$ and define the map $N:K\rightarrow K$ as follows: \[N(u):=(T\circ\rho)(u).\] 
We introduce the following useful notation:
\[(K_i)^{\varphi_i}_{r_i}=\{u\in K_{i}\,:\, \varphi_i(u)< r_{i}\} \quad \text{ and } \quad (\overline{K}_i)^{\varphi_i}_{r_i}=\{u\in K_{i}\,:\, \varphi_i(u)\leq r_{i}\}. \]
Clearly, $\partial\,(K_i)^{\varphi_i}_{r_i}\subset\{u\in K_{i} :  \varphi_i(u)= r_{i}\}$ ($i=1,2$) and, moreover, \[K^{\;\varphi}_{r,R}=\left((K_1)_{R_1}\setminus (\overline{K}_1)^{\varphi_1}_{r_1} \right)\times \left((K_2)_{R_2}\setminus (\overline{K}_2)^{\varphi_2}_{r_2} \right),\] so now the proof follows as those of Theorems \ref{th_CC}, \ref{th_CE} and \ref{th_EE}, replacing $\left\|\cdot\right\|$ with $\varphi_i(\cdot)$ where needed. \qed

\begin{remark}\label{rmk_domain}
	The previous fixed point index computation remains true for an operator $T$ defined in a much more general domain of type $\left(\overline{U}_1\setminus V_1 \right)\times \left(\overline{U}_2\setminus V_2 \right)$, where for each $i\in\{1,2\}$, one has $0\in V_i\subset\overline{V}_i\subset U_i$, $U_i$ and $V_i$ are bounded and relatively open sets in $K_i$ and $\overline{U}_i\setminus V_i$ is a retract of $\overline{U}_i$. Observe that, in particular, $\overline{U}_i\setminus V_i$ is a retract of $\overline{U}_i$ provided that $\partial\,V_i$ is a retract of $\overline{V}_i$, what allows $U_i$ to be an arbitrary bounded open set large enough. It can be immediately deduced that this is the case for $V_i=(K_i)_{r_i}$. Indeed, for a fixed $h_i\in K_i\setminus\{0\}$, one can define the retraction $\rho_i:\overline{V}_i\rightarrow\partial\,V_i$ as
	\[\rho_i(u_i)=r_i\displaystyle\frac{u_i+(r_i-\left\|u_i\right\|)^2 h_i}{\left\|u_i+(r_i-\left\|u_i\right\|)^2 h_i \right\|}. \]
	 Note that in this case the set $U_i$ needs not be the intersection of a ball with the cone $K_i$, which enlarges the applicability of Theorems \ref{th_CC}, \ref{th_CE} and \ref{th_EE}.
	
	In this context, conditions $(i)$ and $(ii)$ above can be written in the following way:
	\begin{enumerate}
		\item[$(i)$] $T_i(u)+\mu\,h_i\neq u_i$ if $u_i\in\partial\,V_i$ and $\mu\geq 0$, and $T_i(u)\neq \lambda\, u_i$ if $u_i\in\partial\,U_i$ and $\lambda\geq 1$;
		\item[$(ii)$] $T_i(u)\neq \lambda\, u_i$ if $u_i\in\partial\,V_i$ and $\lambda\geq 1$, and $T_i(u)+\mu\,h_i\neq u_i$ if $u_i\in\partial\,U_i$ and $\mu\geq 0$; 
	\end{enumerate}
	and they must be satisfied over the set $\left(\overline{U}_1\setminus V_1 \right)\times \left(\overline{U}_2\setminus V_2 \right)$.
\end{remark}

Obviously, from Theorem \ref{th_ind_funct} it follows immediately a new fixed point theorem in the line of {Theorem~\ref{th_KP}}.

\begin{theorem}\label{th_KP_funct}
%
	Assume that there exist continuous concave functionals $\varphi_i:K_i\rightarrow\mathbb{R}_+$ such that $\varphi_i(u)\leq\left\|u\right\|$ for all $u\in K_i$ ($i=1,2$), the set $K^{\;\varphi}_{r,R}$ is nonempty, $T=(T_1,T_2):\overline{K}^{\;\varphi}_{r,R}\rightarrow K$ is a compact map and for each $i\in\{1,2\}$ there exists $h_i\in K_i\setminus\{0\}$ such that one of the following conditions is satisfied in $\overline{K}^{\;\varphi}_{r,R}$:
	\begin{enumerate}
		\item[$(a)$] $T_i(u)+\mu\,h_i\neq u_i$ if $\varphi_i(u_i)=r_i$ and $\mu\geq 0$, and $T_i(u)\neq \lambda\, u_i$ if $\left\|u_i\right\|=R_i$ and $\lambda\geq 1$;
		\item[$(b)$] $T_i(u)\neq \lambda\, u_i$ if $\varphi_i(u_i)=r_i$ and $\lambda\geq 1$, and $T_i(u)+\mu\,h_i\neq u_i$ if $\left\|u_i\right\|=R_i$ and $\mu\geq 0$. 
	\end{enumerate}
	
	Then $T$ has at least a fixed point $u=(u_1,u_2)\in K$ with $r_i< \varphi_i(u_i)$ and $\left\|u_i\right\|< R_i$ $(i=1,2)$.
\end{theorem}

\begin{remark}
	It seems hard to adapt the arguments in the original proof of Theorem \ref{th_KP} due to Precup  \cite{PrecupFPT,PrecupSDC} to demonstrate Theorem \ref{th_KP_funct}. In particular, it may be difficult to reduce the expansive type conditions to compressive ones, since there is no relation between the functional $\varphi_i$ and the corresponding norm. However, our reasonings, based on fixed point index theory, work similarly for both proofs.  
\end{remark}

\subsection{Multiplicity results: a three-solutions type fixed point theorem}

Topological methods have been frequently employed to establish multiplicity results. For instance, we mention the three-solutions fixed point theorems in cones due to Amann \cite{amann} and due to Leggett and Williams \cite{LW}, which are well-known and widely applied in the literature. By using the computations of the fixed point index obtained in the previous section, we can easily derive a multiplicity result in this line.

For simplicity, we will restrict our comments here to compact operators defined in sets of the form $\overline{K}_{r,R}$, given as above in terms of the corresponding norm. Obviously, similar conclusions can be obtained for operators defined in other distinct regions where the previous computations of the fixed point index remain valid, for instance those of type $\overline{K}^{\;\varphi}_{r,R}$, defined by means of appropriated functionals. 

\begin{theorem}\label{th_mult}
	Let $X$ and $Y$ be normed linear spaces, $K_1\subset X$ and $K_2\subset Y$ two cones, $K:=K_1\times K_2$ and $r^{(j)},R^{(j)}\in\R^2_+$, $r^{(j)}=(r^{(j)}_{1},r^{(j)}_{2})$, $R^{(j)}=(R^{(j)}_{1},R^{(j)}_{2})$, with $0<r^{(j)}_i<R^{(j)}_i$ $(i=1,2, \ j=1,2,3)$. Assume that the sets $\overline{K}_{r^{(j)},R^{(j)}}$ are such that 
	\[\overline{K}_{r^{(1)},R^{(1)}}\cup \overline{K}_{r^{(2)},R^{(2)}}\subset \overline{K}_{r^{(3)},R^{(3)}} \quad \text{ and } \quad \overline{K}_{r^{(1)},R^{(1)}}\cap \overline{K}_{r^{(2)},R^{(2)}}=\emptyset. \]
	
	Moreover, assume that $T=(T_1,T_2):\overline{K}_{r^{(3)},R^{(3)}}\rightarrow K$ is a compact map and for each $i\in\{1,2\}$ and each $j\in\{1,2,3\}$ there exists $h_i^j\in K_i\setminus\{0\}$ such that one of the following conditions is satisfied in $\overline{K}_{r^{(j)},R^{(j)}}$:
	\begin{enumerate}
		\item[$(a)$] $T_i(u)+\mu\,h_i^j\neq u_i$ if $\left\|u_i\right\|=r^{(j)}_i$ and $\mu\geq 0$, and $T_i(u)\neq \lambda\, u_i$ if $\left\|u_i\right\|=R^{(j)}_i$ and $\lambda\geq 1$;
		\item[$(b)$] $T_i(u)\neq \lambda\, u_i$ if $\left\|u_i\right\|=r^{(j)}_i$ and $\lambda\geq 1$, and $T_i(u)+\mu\,h_i\neq u_i$ if $\left\|u_i\right\|=R^{(j)}_i$ and $\mu\geq 0$. 
	\end{enumerate}
	Then $T$ has at least three distinct fixed points $\bar{u}^1$, $\bar{u}^2$ and $\bar{u}^3$ such that 
	\[\bar{u}^1\in K_{r^{(1)},R^{(1)}}, \ \ \bar{u}^2\in K_{r^{(2)},R^{(2)}} \ \ \text{ and } \ \ \bar{u}^3\in K_{r^{(3)},R^{(3)}}\setminus\left(\overline{K}_{r^{(1)},R^{(1)}}\cup \overline{K}_{r^{(2)},R^{(2)}} \right). \]
\end{theorem}

\noindent
{\bf Proof.} It follows from Theorems \ref{th_CC}, \ref{th_CE} and \ref{th_EE} that 
\begin{equation}\label{eq_ind_mul}
i_{K}\left(T,K_{r^{(j)},R^{(j)}}\right)=\pm 1 \qquad (j=1,2,3). 
\end{equation}
By the existence property of the fixed point index, one obtains that $T$ has two fixed points $\bar{u}^1\in K_{r^{(1)},R^{(1)}}$ and $\bar{u}^2\in K_{r^{(2)},R^{(2)}}$. Clearly, $\bar{u}^1$ and $\bar{u}^2$ are different, given that $K_{r^{(1)},R^{(1)}}\cap K_{r^{(2)},R^{(2)}}=\emptyset$.

On the other hand, since $\overline{K}_{r^{(1)},R^{(1)}}\cup \overline{K}_{r^{(2)},R^{(2)}}\subset \overline{K}_{r^{(3)},R^{(3)}}$ and $\overline{K}_{r^{(1)},R^{(1)}}\cap \overline{K}_{r^{(2)},R^{(2)}}=\emptyset$, one has
\[K_{r^{(3)},R^{(3)}}\setminus\left(\overline{K}_{r^{(1)},R^{(1)}}\cup \overline{K}_{r^{(2)},R^{(2)}} \right)\neq\emptyset. \]
Hence, the additivity property of the fixed point index combined with \eqref{eq_ind_mul} imply that
\[i_{K}\left(T,K_{r^{(3)},R^{(3)}}\setminus\left(\overline{K}_{r^{(1)},R^{(1)}}\cup \overline{K}_{r^{(2)},R^{(2)}} \right)\right)=i_{K}\left(T,K_{r^{(3)},R^{(3)}}\right)-i_{K}\left(T,K_{r^{(1)},R^{(1)}}\right)-i_{K}\left(T,K_{r^{(2)},R^{(2)}}\right) \]
is an odd number and thus, 
\[i_{K}\left(T,K_{r^{(3)},R^{(3)}}\setminus\left(\overline{K}_{r^{(1)},R^{(1)}}\cup \overline{K}_{r^{(2)},R^{(2)}} \right)\right)\neq 0. \]
Therefore, the operator $T$ has a third fixed point $\bar{u}^3\in K_{r^{(3)},R^{(3)}}\setminus\left(\overline{K}_{r^{(1)},R^{(1)}}\cup \overline{K}_{r^{(2)},R^{(2)}} \right)$.
\qed

\subsection{Comparison with classical results: Poincar\'e-Miranda theorem} 

In this section we show that the finite-dimensional version of Krasnosel'ski\u{\i}-Precup fixed point theorem is equivalent to the classical Poincar\'e-Miranda theorem concerning the existence of a zero of a nonlinear mapping in a rectangle of $\R^n$. Poincar\'e-Miranda theorem is the $n$-dimensional version of Bolzano intermediate value theorem and it is well-known its equivalence to Brouwer fixed point theorem (see \cite{Maw} for historical details). It has received a renewed interest in recent years, see for instance \cite{fel,FG,Frank,Maw_CE,S} and the references therein for some new versions, generalizations and applications of this result.

First of all, let us rewrite Theorem \ref{th_KPn} in the particular case in which $X=\R$ with the usual norm $\left\|\cdot\right\|=\left|\cdot\right|$, $K_i=\R_{+}$ (and thus $K=\R_+^n$) and $h_i=1$ for every $i\in\{1,2,\dots,n \}$. 

\begin{theorem}\label{th_KP_Rn}
	Let $r,R\in\R^n_+$, $r=(r_{1},\dots,r_{n})$, $R=(R_{1},\dots,R_{n})$, with $0<r_i<R_i$ $(i=1,\dots,n)$.
	
	Assume that $f=(f_1,\dots,f_n):\overline{K}_{r,R}\rightarrow \R_+^n$ is a continuous function and for each $i\in\{1,\dots,n\}$ one of the following conditions is satisfied in $\overline{K}_{r,R}$:
	\begin{enumerate}
		\item[$(a)$] $f_i(x)\geq x_i$ if $\left|x_i\right|=r_i$, and $f_i(x)\leq x_i$ if $\left|x_i\right|=R_i$;
		\item[$(b)$] $f_i(x)\leq x_i$ if $\left|x_i\right|=r_i$, and $f_i(x)\geq x_i$ if $\left|x_i\right|=R_i$. 
	\end{enumerate}
	Then $f$ has at least a fixed point $\bar{x}=(\bar{x}_1,\dots,\bar{x}_n)\in \R_+^n$ with $r_i\leq\left|\bar{x}_i\right|\leq R_i$ $(i=1,\dots,n)$.
\end{theorem}

\begin{remark}
	In this setting, the set $\overline{K}_{r,R}$ is the following rectangle in $\R_+^n$: $\mathcal{R}=[r_1,R_1]\times\cdots\times[r_n,R_n]$. Furthermore, the points $x=(x_1,\dots,x_n)\in \overline{K}_{r,R}$ satisfying that $\left|x_i\right|=r_i$ (respectively, $\left|x_i\right|=R_i$) are in fact those with $x_i=r_i$ (resp. $x_i=R_i$).	
\end{remark}

Let us now recall the Poincar\'e-Miranda theorem. We shall see that it can be directly deduced from Theorem \ref{th_KP_Rn}.

\begin{theorem}[Poincar\'e-Miranda]
	Let $\mathcal{R}=[a_1,b_1]\times\cdots\times[a_n,b_n]$ be a rectangle in $\R^n$. Assume that $g=(g_1,\dots,g_n):\mathcal{R}\rightarrow\R^n$ is a continuous function and for each $i\in\{1,\dots,n\}$ one of the following conditions is satisfied in $\mathcal{R}$:
	\begin{enumerate}
		\item[$(a)$] $g_i(x)\geq 0$ if $x_i=a_i$, and $g_i(x)\leq 0$ if $x_i=b_i$;
		\item[$(b)$] $g_i(x)\leq 0$ if $x_i=a_i$, and $g_i(x)\geq 0$ if $x_i=b_i$. 
	\end{enumerate}
	Then there exists $\bar{x}\in\mathcal{R}$ such that $g(\bar{x})=0$. 
\end{theorem}

\noindent
{\bf Proof.} Let us divide the proof in two cases.

\textit{Case 1: $a_i>0$ for every $i\in\{1,2,\dots,n \}$}. Since $g$ is a continuous function in the compact set $\mathcal{R}$, it is bounded. Then choose $\lambda\in(0,1)$ such that $\lambda\left|g_i(x)\right|\leq a$ for all $x\in \mathcal{R}$ and all $i\in\{1,\dots,n \}$, with $a:=\min\{a_1,a_2,\dots,a_n\}$. Now define the function $f:\mathcal{R}\rightarrow\R^n$ as follows
\[f(x)=\lambda\, g(x)+x. \]

Note that $f(\mathcal{R})\subset\R_+^n$. Indeed, for a fixed $i\in\{1,\dots,n \}$, one has that $\lambda\,g_i(x)\geq -a\geq-a_i$ for all $x\in \mathcal{R}$ and thus 
$f_i(x)=\lambda\,g_i(x)+x_i\geq -a_i+x_i\geq 0$ for all $x\in\mathcal{R}$. Moreover, condition $(a)$ implies that 
\begin{enumerate}
	\item[$(a^*)$] $f_i(x)\geq x_i$ if $x_i=a_i$, and $f_i(x)\leq x_i$ if $x_i=b_i$,
\end{enumerate}
and, similarly, it follows from $(b)$ that 
\begin{enumerate}
	\item[$(b^*)$] $f_i(x)\leq x_i$ if $x_i=a_i$, and $f_i(x)\geq x_i$ if $x_i=b_i$.
\end{enumerate}
Hence, Theorem \ref{th_KP_Rn} ensures that $f$ has a fixed point $\bar{x}\in \mathcal{R}$, that is, \[\bar{x}=f(\bar{x})=\lambda\,g(\bar{x})+\bar{x},\] and so $\bar{x}$ is a zero of $g$. 

\textit{Case 2: $a_i\leq 0$ for some $i\in\{1,2,\dots,n \}$}. A trivial translation moves the rectangle $\mathcal{R}$ to the interior of $\R_+^n$ and the conclusion follows then from Case 1. 
\qed

\begin{remark}
	Clearly, Theorem \ref{th_KP_Rn} can be seen as a consequence of Poincar\'e-Miranda existence theorem. Indeed, assume that $f:\mathcal{R}\rightarrow \R_+^n$ is a continuous function under the assumptions of Theorem \ref{th_KP_Rn} and let $g:\mathcal{R}\rightarrow \R^n$ be the continuous map defined by $g(x)=x-f(x)$, $x\in\mathcal{R}$. It follows from Poincar\'e-Miranda theorem that $g$ has a zero in $\mathcal{R}$, which obviously is a fixed point of $f$. 
\end{remark}

\section{Applications}

\subsection{Hammerstein systems}

Consider the following system of Hammerstein type equations
\begin{equation}\label{eq_Ham}
	\begin{array}{r}
	u_1(t)=\displaystyle\int_{0}^{1} k_1(t,s)g_1(s)f_1(u_1(s),u_2(s))\,ds, \\[0.3cm]
	u_2(t)=\displaystyle\int_{0}^{1} k_2(t,s)g_2(s)f_2(u_1(s),u_2(s))\,ds,
	\end{array}
\end{equation}
where $I:=[0,1]$ and for each $i\in\{1,2\}$ the following assumptions are satisfied:
\begin{enumerate}
	\item[$(H_1)$] the kernel $k_i:I^2\rightarrow\mathbb{R}_+$ is continuous;
	\item[$(H_2)$] the function $g_i:I\rightarrow \mathbb{R}_+$ is measurable;
	\item[$(H_3)$] there exist an interval $[a,b]\subset I$ and a function $\Phi_i:I\rightarrow\mathbb{R}_+$ such that
	\[\Phi_i\,g_i\in L^1(I), \qquad \int_{a}^{b}\Phi_i(s)g_i(s)\,ds>0, \]
	and a constant $c_i\in (0,1]$ satisfying
	\[\begin{array}{rll}
	k_i(t,s)&\leq \Phi_i(s)  & \quad \text{for all } t,s\in I, \\ c_i\,\Phi_i(s)&\leq k_i(t,s) & \quad \text{for all } t\in[a,b], \ s\in I;
	\end{array}\]
	\item[$(H_4)$] the function $f_i:\mathbb{R}_+^2\rightarrow\mathbb{R}_+$ is continuous.
\end{enumerate}	
 
Let us consider the Banach space of continuous functions $X=\mathcal{C}(I)$ endowed with the usual maximum norm $\left\|v\right\|_{\infty} =\max_{t\in I}\left|v(t)\right|$ and the cones 
\[K_i=\left\{v\in \mathcal{C}(I)\,:\,v(t)\geq 0 \text{ for all } t\in I \text{ and } \min_{t\in[a,b]}v(t)\geq c_i\left\|v\right\|_{\infty} \right\} \qquad (i=1,2). \] 

In order to prove the existence of positive solutions of the system of integral equations \eqref{eq_Ham}, we look for fixed points of the operator $T:K\rightarrow K$, $T=(T_1,T_2)$, defined as 
\begin{equation}\label{eqT_Ham}
T_i(u_1,u_2)(t):= \displaystyle\int_{0}^{1} k_i(t,s)g_i(s)f_i(u_1(s),u_2(s))\,ds \qquad (i=1,2),
\end{equation}
where $K:=K_1\times K_2$ is a cone in $X^2$. Under assumptions $(H_1)$--$(H_4)$, it can be proven by means of standard arguments (see, for instance, \cite{fig_tojo,Inf}) that $T$ maps the cone $K$ into itself and it is completely continuous, i.e., $T$ is continuous and maps bounded sets into relatively compact ones.

Now, let us fix some notations:
\[A_i:=\inf_{t\in[a,b]}\int_{a}^{b}k_i(t,s)g_i(s)\,ds, \quad B_i:=\sup_{t\in I}\int_{0}^{1}k_i(t,s)g_i(s)\,ds \quad (i=1,2). \]
In addition, for $\alpha_i,\beta_i>0$, $\alpha_i\neq \beta_i$, $i=1,2$, denote
\begin{align*}
m_1^{\alpha,\beta}&:=\min\{f_1(u_1,u_2)\,:\,c_1\,\beta_1\leq u_1\leq \beta_1, \ c_2\,r_2\leq u_2\leq R_2 \}, \\
m_2^{\alpha,\beta}&:=\min\{f_2(u_1,u_2)\,:\,c_1\,r_1\leq u_1\leq R_1, \ c_2\,\beta_2\leq u_2\leq \beta_2 \}, \\
M_1^{\alpha,\beta}&:=\max\{f_1(u_1,u_2)\,:\,0\leq u_1\leq \alpha_1, \ 0\leq u_2\leq R_2 \}, \\
M_2^{\alpha,\beta}&:=\max\{f_2(u_1,u_2)\,:\,0\leq u_1\leq R_1, \ 0\leq u_2\leq \alpha_2 \}, 
\end{align*}
where $r_i:=\min\{\alpha_i,\beta_i \}$ and $R_i:=\max\{\alpha_i,\beta_i \}$.

We are in a position to establish an existence result for the system of Hammerstein type equations \eqref{eq_Ham} as a consequence of Theorem \ref{th_KP_index}.

\begin{theorem}\label{th_Ham}
	Assume that conditions $(H_1)$--$(H_4)$ are fulfilled. Moreover, suppose that there exist positive numbers $\alpha_i,\beta_i>0$ with $\alpha_i\neq \beta_i$, $i=1,2$, such that
	\begin{equation}\label{eq_CE_Ham}
	A_i\,m_i^{\alpha,\beta}>\beta_i, \qquad B_i\,M_i^{\alpha,\beta}<\alpha_i \quad (i=1,2).
	\end{equation}	
	Then the system \eqref{eq_Ham} has at least one positive solution $(u_1,u_2)\in K$ such that $r_i<\left\|u_i\right\|_{\infty}<R_i$ $(i=1,2)$.
\end{theorem}

\noindent
{\bf Proof.} Consider the operator $T:\overline{K}_{r,R}\rightarrow K$, $T=(T_1,T_2)$, defined as in \eqref{eqT_Ham}. Let us check that for every $u=(u_1,u_2)\in \overline{K}_{r,R}$ and $i\in\{1,2\}$ the following conditions are satisfied:
\begin{enumerate}
	\item[$1)$] $T_i(u)+\mu\,{\pmb 1}\neq u_i$ if $\left\|u_i\right\|_{\infty}=\beta_i$ and $\mu\geq0$ (where ${\pmb 1}$ denotes the constant function equal to one); 
	\item[$2)$] $T_i(u)\neq \lambda\, u_i$ if $\left\|u_i\right\|_{\infty}=\alpha_i$ and $\lambda\geq 1$.
\end{enumerate}

First, to prove $1)$, assume to the contrary that there exist $u=(u_1,u_2)\in \overline{K}_{r,R}$ with $\left\|u_i\right\|_{\infty}=\beta_i$ and $\mu\geq0$ such that $T_i(u)+\mu\,{\pmb 1}= u_i$. Then we have 
\[u_i(t)=\displaystyle\int_{0}^{1} k_i(t,s)g_i(s)f_i(u(s))\,ds+\mu. \]
Since $u\in \overline{K}_{r,R}$ with $\left\|u_i\right\|_{\infty}=\beta_i$, it follows from the definition of the cone $K$ that $c_i\,\beta_i\leq u_i(t)\leq \beta_i$ and $c_j\,r_j\leq u_j(t)\leq R_j$ ($j\neq i$) for all $t\in[a,b]$. Hence, for $t\in[a,b]$, we deduce from \eqref{eq_CE_Ham} that
\[u_i(t)\geq \int_{a}^{b} k_i(t,s)g_i(s)f_i(u(s))\,ds\geq m_i^{\alpha,\beta}\int_{a}^{b} k_i(t,s)g_i(s)\,ds\geq A_i\,m_i^{\alpha,\beta}>\beta_i, \]
a contradiction.

Now, let us show that $\left\|T_i(u)\right\|_{\infty}<\alpha_i$ for all $u\in \overline{K}_{r,R}$ with $\left\|u_i\right\|_{\infty}=\alpha_i$, which clearly implies property $2)$. Note that, if $u\in \overline{K}_{r,R}$ with $\left\|u_i\right\|_{\infty}=\alpha_i$, then $0\leq u_i(t)\leq \alpha_i$ and $0\leq u_j(t)\leq R_j$ ($j\neq i$) for all $t\in I$. Thus, we get, for $t\in I$,
\[T_i(u)(t)=\displaystyle\int_{0}^{1} k_i(t,s)g_i(s)f_i(u(s))\,ds\leq M_i^{\alpha,\beta}\int_{0}^{1} k_i(t,s)g_i(s)\,ds\leq B_i\,M_i^{\alpha,\beta}. \]
Taking the maximum over $I$ and applying condition \eqref{eq_CE_Ham}, we obtain $\left\|T_i(u)\right\|_{\infty}\leq B_i\,M_i^{\alpha,\beta}<\alpha_i$.

Finally, the result follows from Theorem \ref{th_KP_index} with $r_i=\min\{\alpha_i,\beta_i \}$ and $R_i=\max\{\alpha_i,\beta_i \}$.
\qed

\begin{remark}
	When we have additional information concerning the monotonicity of the nonlinearities $f_1$ and $f_2$, we may rewrite condition \eqref{eq_CE_Ham} in terms of the values of $f_1$ and $f_2$ at certain points. For instance, if $f_1$ and $f_2$ are nondecreasing on $[0,R_1]\times[0,R_2]$, then \eqref{eq_CE_Ham} is equivalent to
	\begin{align*}
	A_1\,f_1(c_1\,\beta_1,c_2\,r_2)&>\beta_1, \qquad B_1\,f_1(\alpha_1,R_2)<\alpha_1, \\ A_2\,f_2(c_1\,r_1,c_2\,\beta_2)&>\beta_2, \qquad B_2\,f_2(R_1,\alpha_2)<\alpha_2.
	\end{align*}
\end{remark}

Now we will apply Theorem \ref{th_mult} in order to give sufficient conditions for the existence of at least three positive solutions for the system \eqref{eq_Ham}.

\begin{theorem}\label{th_mult_Ham}
	Assume that conditions $(H_1)$--$(H_4)$ hold. Moreover, suppose that there exist $\alpha_i^j,\beta_i^j>0$ with $\alpha_i^j\neq \beta_i^j$, $i=1,2$, $j=1,2,3$, such that 
	$\alpha_i^1,\alpha_i^2,\beta_i^1,\beta_i^2\in[r_i^3,R_i^3]$ for $i\in\{1,2\}$,  there exists $i\in\{1,2\}$ such that $R_i^1<r_i^2$,  
	where $r_i^j:=\min\{\alpha_i^j,\beta_i^j \}$ and $R_i^j:=\max\{\alpha_i^j,\beta_i^j \}$,	and
	\begin{equation}\label{eq_CE_Ham_mult}
	A_i\,m_i^{\alpha^j,\beta^j}>\beta_i^j, \qquad B_i\,M_i^{\alpha^j,\beta^j}<\alpha_i^j \quad (i=1,2, \ j=1,2,3).
	\end{equation}	
	Then the system \eqref{eq_Ham} has at least three positive solutions.
\end{theorem}

\noindent
{\bf Proof.} Consider the operator $T:\overline{K}_{r^3,R^3}\rightarrow K$, $T=(T_1,T_2)$, defined as in \eqref{eqT_Ham}. 

Since $\alpha_i^1,\alpha_i^2,\beta_i^1,\beta_i^2\in[r_i^3,R_i^3]$ for $i\in\{1,2\}$, we have that $\overline{K}_{r^{1},R^{1}}\cup \overline{K}_{r^{2},R^{2}}\subset \overline{K}_{r^{3},R^{3}}$. In addition, the existence of some $i\in\{1,2\}$ such that $R_i^1<r_i^2$ clearly implies that $\overline{K}_{r^{1},R^{1}}\cap \overline{K}_{r^{2},R^{2}}=\emptyset$.

Finally, conditions $(a)$ and $(b)$ in Theorem \ref{th_mult} can be verified by using the inequalities given in \eqref{eq_CE_Ham_mult} and completely analogous reasonings to those in the proof of Theorem \ref{th_Ham}. Therefore, we reach the thesis as a consequence of Theorem \ref{th_mult}.
\qed


\begin{example}
Consider the following system of second-order equations with Dirichlet boundary conditions	
\begin{equation}\label{ex_2or}
\begin{array}{rl}
-u_1''&=f_1(u_1,u_2), \\ -u_2''&=f_2(u_1,u_2), \\ u_1(0)&=u_1(1)=0=u_2(0)=u_2(1),
\end{array}
\end{equation}
with $f_1(u_1,u_2)=h(u_1)(1+\sin^2(u_2))$, $f_2(u_1,u_2)=u_2^2(1+\sin^2(u_1))$ and 
\[h(u_1)=\left\{\begin{array}{ll}
\sqrt[3]{u_1}, & \quad \text{if } u_1\in[0,1], \\ {u_1}^3, & \quad \text{if } u_1\in(1,10), \\ \sqrt[3]{u_1-10}+1000, & \quad \text{if } u_1\in[10,+\infty).
\end{array} \right. \]

We can associate to \eqref{ex_2or} a system of Hammerstein type equations of the form \eqref{eq_Ham} where the kernels are given by the corresponding Green's function
\[k_1(t,s)=k_2(t,s)=\left\{\begin{array}{ll}
s(1-t), & \quad \text{if } s\leq t, \\ t(1-s), & \quad \text{if } s>t,
\end{array} \right. \]
and $g_1(t)=g_2(t)=1$ for all $t\in I$. It is well-known (see, for instance, \cite{Inf}) that condition $(H_3)$ holds if we take
\[\Phi_i(s)=s(1-s), \quad [a,b]=[1/4,3/4], \quad c_i=1/4 \quad (i=1,2). \]
This choice leads to $A_i=1/16$ and $B_i=1/8$ for $i=1,2$.

Now, it is a routine to check that condition \eqref{eq_CE_Ham_mult} is satisfied with $\beta_1^1=\beta_1^3=2^{-9}$, $\alpha_1^1=2^{-2}$, $\beta_1^2=2^6$, $\alpha_1^2=\alpha_1^3=2^9+10$, $\alpha_2^j=2$ and $\beta_2^j=2^9$ for $j=1,2,3$. Therefore, Theorem \ref{th_mult_Ham} guarantees that the system \eqref{ex_2or} has at least three positive solutions $(u_1,u_2)$, $(v_1,v_2)$ and $(w_1,w_2)$ with the following localizations
\begin{align*}
\frac{1}{512}&<\left\|u_1\right\|_{\infty}<\frac{1}{4}, \ \ 2<\left\|u_2\right\|_{\infty}<512, \\ 64&<\left\|v_1\right\|_{\infty}<522, \ \ 2<\left\|v_2\right\|_{\infty}<512, \\ \frac{1}{4}&\leq\left\|w_1\right\|_{\infty}\leq 64, \ \ 2<\left\|w_2\right\|_{\infty}<512.
\end{align*}
We emphasize that the second component of all the three solutions is situated in the same region, so the multiplicity is obtained due to we are able to localize their first component in distinct sets.
\end{example}

\subsection{Radial solutions of $(p_1,p_2)$-Laplacian systems}

In this section, we consider the existence of positive radial solutions for the $(p_1,p_2)$-Laplacian system 
\begin{equation}\label{eq_pLap}
\begin{array}{rll}
-\Delta_{p_1} u&=f_1(u,v) & \quad \text{ in } B, \\ -\Delta_{p_2} v&=f_2(u,v) & \quad \text{ in } B, \\ u&=0=v & \quad \text{ on } \partial B,
\end{array}
\end{equation}
where $\Delta_p u={\rm div}\left(\left\|\nabla u\right\|^{p-2}\nabla u\right)$, $B$ is the unit open ball in $\mathbb{R}^n$ centered at origin, $p_1,p_2>n\geq 2$ and $f_1,f_2:\mathbb{R}_+^2\rightarrow\mathbb{R}_+$ are continuous and nondecreasing functions (that is, if $(u_1,v_1),(u_2,v_2)\in \mathbb{R}_+^2$ with $u_1\leq u_2$ and $v_1\leq v_2$, then $f_i(u_1,v_1)\leq f_i(u_2,v_2)$ for $i=1,2$).

Setting, as usual, $r=\left\|x\right\|$, $u(x)=u_1(r)$ and $v(x)=u_2(r)$, the Dirichlet system \eqref{eq_pLap} is reduced to the following system of ordinary differential equations with mixed boundary conditions
\begin{equation}\label{eq_rad}
\begin{array}{rl}
-\left[r^{n-1}\phi_{p_1}(u_1') \right]'&=r^{n-1}f_1(u_1,u_2)  \quad \text{ in } (0,1), \\ -\left[r^{n-1}\phi_{p_2}(u_2') \right]'&=r^{n-1}f_2(u_1,u_2)  \quad \text{ in } (0,1), \\ u_1'(0)&=u_1(1)=0=u_2'(0)=u_2(1), 
\end{array}
\end{equation}
where $\phi_p(t):=\left|t\right|^{p-2}t$ is the $p$-Laplacian homeomorphism. We will look for positive solutions of \eqref{eq_rad}, that is, radially symmetric solutions of \eqref{eq_pLap}.

A Harnack type inequality has been established in \cite{PPV} for problem 
\[-\left[r^{n-1}\phi_p(v') \right]'=r^{n-1}h(r,v)  \quad \text{ in } (0,1), \quad v'(0)=v(1)=0, \]
in terms of the \textit{energetic} norm. By using H\"{o}lder inequality, one can derive a Harnack type inequality in terms of the usual max-norm, see \cite{HOrP}. The result can be summarized as follows.

\begin{lemma}\label{lem_H}
	Let $p>n$. Every function $v\in C^{1}\left[ 0,1\right] $ with $r^{n-1}\phi_p(v')\in C^{1}\left[ 0,1\right] $ and $\left[r^{n-1}\phi_p(v') \right]' \leq 0$ on $\left[
	0,1\right]$ satisfies that $v^{\prime }\leq 0$ on $\left[ 0,1\right] .$
	If, in addition, $-r^{1-n}\left[r^{n-1}\phi_p(v') \right]'$ is nonincreasing on $(0,1],$ then
	\begin{equation*}
	v\left( r\right) \geq  \frac{p-n}{p-1}\left(
	1-r\right) r^{\frac{n}{p-1}}\left\Vert v\right\Vert _{\infty},\ \ \  r\in %
	\left[ 0,1\right]. 
	\end{equation*}%
\end{lemma}

Let us consider the following cones in the space of continuous functions $\mathcal{C}(I)$, with $I:=[0,1]$,
\[K_i=\left\{ v\in \mathcal{C}(I)\,:\,v\geq 0 \text{ on } I,\ v \text{ is nonincreasing on } I \text{ and } \min_{r\in[a,b]}v(r)\geq c_i\left\|v\right\|_{\infty} \right\} \qquad (i=1,2), \]
where $[a,b]\subset(0,1)$ and $c_i:=\frac{p_i-n}{p_i-1}\left(1-b\right) a^{\frac{n}{p_i-1}}$, $i=1,2$. As before, we define the cone $K:=K_1\times K_2$ in the product space.

We will look for positive solutions of problem \eqref{eq_rad} as fixed points of the operator $T=(T_1,T_2):K\rightarrow K$ defined as
\begin{equation}\label{eq_Trad}
T_i(u_1,u_2)(r)=\int_{r}^{1}\phi_{p_i}^{-1}\left(\frac{1}{s^{n-1}}\int_{0}^{s}\tau^{n-1}f_i(u_1(\tau),u_2(\tau))\,d\tau\right)ds, \qquad (i=1,2).
\end{equation}

The operator $T$ is well-defined, that is, it maps the cone $K$ into itself. Indeed, take $u=(u_1,u_2)\in K$ and let us show that $v_i:=T_i(u)\in K_i$, for $i=1,2$. Clearly, $v_i\in \mathcal{C}(I)$ and $v_i\geq 0$ on $I$, since $f_i$ is continuous and nonnegative. Moreover, $\left[r^{n-1}\phi_{p_i}(v'_i) \right]' \leq 0$ on $I$, so $v_i$ is nonincreasing on $I$, as a consequence of Lemma \ref{lem_H}. On the other hand, $u_1$ and $u_2$ are nonincreasing and $f_i$ is nondecreasing, which implies that the map $r\mapsto f_i(u_1(r),u_2(r))$ is nonincreasing on $I$ and then so is $-r^{1-n}\left[r^{n-1}\phi_{p_i}(v'_i) \right]'=f_i(u_1(r),u_2(r))$. Hence, again by Lemma \ref{lem_H}, one has that $\min_{r\in[a,b]}v_i(r)\geq c_i\left\|v_i\right\|_{\infty}$. In conclusion, $v_i\in K_i$, as desired.

It is a routine to check that $T$ is completely continuous.

Now, let us define a continuous concave functional on $K_i$,  $\varphi_i:K_i\rightarrow\mathbb{R}_+$, as follows
\[\varphi_i(v)=\min_{r\in[a,b]}v(r), \qquad i=1,2. \]
We intend to apply Theorem \ref{th_KP_funct} in order to obtain sufficient conditions for the existence of fixed points of $T$ located in a set of the form $K^{\;\varphi}_{r,R}$. Note that this set can also be written as $K^{\;\varphi}_{r,R}=(U_1\setminus\overline{V}_1)\times(U_2\setminus\overline{V}_2)$, where
\[V_i=\left\{u\in K_i:\min_{r\in[a,b]}u(r)<r_i \right\}, \qquad U_i=\{u\in K_i:\left\|u\right\|_{\infty}<R_i \} \qquad (i=1,2). \]
The bounded open sets $V_i$ were introduced by Lan in \cite{Lan} and later employed by several authors, see \cite{Inf} and the references therein.

\begin{theorem}\label{th_rad}
	Assume that there exist $\alpha_i,\beta_i>0$ with $\beta_i/c_i<\alpha_i$, $i=1,2$, such that
	\begin{equation}\label{eq_C_rad}
	f_i(\beta_1,\beta_2)>\frac{\beta_i^{p_i-1}}{(b-a)a^{n-1}(1-b)^{p_i-1}}, \qquad f_i(\alpha_1,\alpha_2)<\alpha_i^{p_i-1} \qquad (i=1,2).
	\end{equation}
	Then the system \eqref{eq_rad} has at least one positive solution $(u_1,u_2)\in K$ such that
	\[\beta_i<\varphi_i(u_i) \quad \text{and} \quad \left\|u_i\right\|_{\infty}<\alpha_i, \quad i=1,2. \]
\end{theorem}

\noindent
{\bf Proof.} Consider the operator $T=(T_1,T_2):\overline{K}_{r,R}^{\varphi}\rightarrow K$ defined as in \eqref{eq_Trad}, with $r_i=\beta_i$ and $R_i=\alpha_i$, $i=1,2$. Let us check that it fulfills the assumptions of Theorem \ref{th_KP_funct}.

First, fix $i\in\{1,2\}$ and take $u=(u_1,u_2)\in  \overline{K}_{r,R}^{\varphi}$ with $\varphi_i(u_i)=r_i$. Then $(u_1(r),u_2(r))\geq (r_1,r_2)$ for all $r\in [a,b]$ and thus, by the monotonicity assumption on $f_i$, we have that $f_i(u_1(r),u_2(r))\geq f_i(r_1,r_2)$ for all $r\in[a,b]$. Hence, for $r\in[a,b]$,
\begin{align*}
T_i(u)(r)&\geq \int_{b}^{1}\phi_{p_i}^{-1}\left(\frac{1}{s^{n-1}}\int_{0}^{s}\tau^{n-1}f_i(u_1(\tau),u_2(\tau))\,d\tau\right)ds \\ &\geq\int_{b}^{1}\phi_{p_i}^{-1}\left(\frac{1}{s^{n-1}}\int_{a}^{b}\tau^{n-1}f_i(u_1(\tau),u_2(\tau))\,d\tau\right)ds \\ &\geq\int_{b}^{1}\phi_{p_i}^{-1}\left(\frac{1}{s^{n-1}}\int_{a}^{b}\tau^{n-1}f_i(r_1,r_2)\,d\tau\right)ds \\ &\geq (1-b)\phi_{p_i}^{-1}\left((b-a)a^{n-1}f_i(r_1,r_2) \right)>r_i,
\end{align*}
where the last inequality follows from \eqref{eq_C_rad}. This clearly implies that $T_i(u)+\mu\,{\pmb 1}\neq u_i$ if $u\in  \overline{K}_{r,R}^{\varphi}$ with $\varphi_i(u_i)=r_i$ and $\mu\geq 0$.

Suppose now that $u\in  \overline{K}_{r,R}^{\varphi}$ with $\left\|u_i\right\|=R_i$ for some $i\in\{1,2\}$. Then $f_i(u_1(r),u_2(r))\leq f_i(R_1,R_2)$ for every $r\in I$ and so, by \eqref{eq_C_rad}, we have
\[\left\|T_i(u)\right\|_{\infty}\leq \int_{0}^{1}\phi_{p_i}^{-1}\left(\frac{1}{s^{n-1}}\int_{0}^{s}\tau^{n-1}f_i(u_1(\tau),u_2(\tau))\,d\tau\right)ds\leq \phi_{p_i}^{-1}\left(f_i(R_1,R_2) \right)<R_1. \]
Therefore, alternative $(a)$ in Theorem \ref{th_KP_funct} holds and hence we reach the thesis.
\qed

\begin{remark}
	Assumption \eqref{eq_C_rad} is guaranteed by the following asymptotic conditions: for every $i\in\{1,2\}$, 
	\[\lim_{u_i\to 0}\frac{f_i(u_1,u_2)}{u_i^{p_i-1}}=+\infty \quad \text{and} \quad \lim_{u_i\to \infty}\frac{f_i(u_1,u_2)}{u_i^{p_i-1}}=0 \]
	uniformly with respect to $u_j$, $j\neq i$.
	
	In this case, it is said that both functions $f_1$ and $f_2$ are superlinear at $0$ and sublinear at infinity with respect to $\phi_{p_1}$ and $\phi_{p_2}$, respectively.
\end{remark}

\begin{remark}
	Under the assumptions of Theorem \ref{th_rad}, both operators $T_1$ and $T_2$ are compressive. Notice that the behaviors compressive-expansive and expansive-expansive are also possible:
	\begin{enumerate}
		\item (Compressive-expansive) Assume that there exist $\alpha_i,\beta_i>0$, $i=1,2$, with $\beta_1/c_1<\alpha_1$ and $\alpha_2<\beta_2$, such that
		\begin{equation*}
		\begin{array}{ll}
		f_1(\beta_1,c_2\,\alpha_2)>\displaystyle\frac{\beta_1^{p_1-1}}{(b-a)a^{n-1}(1-b)^{p_1-1}}, & \qquad f_1(\alpha_1,\beta_2/c_2)<\alpha_1^{p_1-1}, \\ f_2(\beta_1,\beta_2)>\displaystyle\frac{\beta_2^{p_2-1}}{(b-a)a^{n-1}(1-b)^{p_2-1}}, & \qquad f_2(\alpha_1,\alpha_2)<\alpha_2^{p_2-1}.
		\end{array}
		\end{equation*}
		Then the operator $T=(T_1,T_2)$ defined as in \eqref{eq_Trad} has at least one fixed point in $\left(\overline{U}_1\setminus V_1 \right)\times \left(\overline{U}_2\setminus V_2 \right)$, where
		\begin{equation*}
		\begin{array}{ll}
			V_1=\{u\in K_1:\varphi_1(u)<\beta_1 \}, & \qquad U_1=\{u\in K_1:\left\|u\right\|_{\infty}<\alpha_1 \}, \\ 
			V_2=\left\{u\in K_2:\left\|u\right\|_{\infty}<\alpha_2 \right\}, & \qquad U_2=\{u\in K_2:\varphi_2(u)<\beta_2 \}.
		\end{array}
		\end{equation*}
		In this case, the operator $T_1$ is compressive and $T_2$ is expansive on $\left(\overline{U}_1\setminus V_1 \right)\times \left(\overline{U}_2\setminus V_2 \right)$.
		\item (Expansive-expansive) Assume that there exist $\alpha_i,\beta_i>0$, with $\alpha_i<\beta_i$, $i=1,2$, such that 
		\begin{equation*}
		\begin{array}{ll}
		f_1(\beta_1,c_2\,\alpha_2)>\displaystyle\frac{\beta_1^{p_1-1}}{(b-a)a^{n-1}(1-b)^{p_1-1}}, & \qquad f_1(\alpha_1,\beta_2/c_2)<\alpha_1^{p_1-1}, \\ f_2(c_1\,\alpha_1,\beta_2)>\displaystyle\frac{\beta_2^{p_2-1}}{(b-a)a^{n-1}(1-b)^{p_2-1}}, & \qquad f_2(\beta_1/c_1,\alpha_2)<\alpha_2^{p_2-1}.
		\end{array}
		\end{equation*}
		Then the operator $T=(T_1,T_2)$ has at least one fixed point in $\left(\overline{U}_1\setminus V_1 \right)\times \left(\overline{U}_2\setminus V_2 \right)$, where
		\[V_i=\left\{u\in K_i:\left\|u\right\|_{\infty}<\alpha_i \right\}, \qquad U_i=\{u\in K_i:\varphi_i(u)<\beta_i \} \qquad (i=1,2). \]
		Note that both operators $T_1$ and $T_2$ are of expansive type on $\left(\overline{U}_1\setminus V_1 \right)\times \left(\overline{U}_2\setminus V_2 \right)$.
	\end{enumerate}
The conclusion can be obtained essentially as in the proof of Theorem \ref{th_rad}, taking into account Remark \ref{rmk_domain}. 
\end{remark}

\section*{Acknowledgements}

Jorge Rodr\'iguez--L\'opez was partially supported by Xunta de Galicia (Spain), project ED431C 2019/02 and AIE, Spain and FEDER, grant PID2020-113275GB-I00.

\end{document}